\documentclass[11pt]{article}
\usepackage{latexsym,amssymb}
\usepackage{psfrag}
\usepackage{amsmath}

\newtheorem{Theorem}{\bf Theorem}[section]
\newtheorem{Lemma}{\bf Lemma}[section]
\newtheorem{Proposition}{\bf Proposition}[section]
\newtheorem{Corollary}{\bf Corollary}[section]
\newtheorem{Remark}{\bf Remark}[section]
\newtheorem{Example}{\bf Example}[section]
\newtheorem{Definition}{\bf Definition}[section]

\newenvironment{theorem}{\begin{Theorem}$\!\!\!$}{\end{Theorem}}
\newenvironment{lemma}{\begin{Lemma}$\!\!\!$}{\end{Lemma}}

\newenvironment{corollary}{\begin{Corollary}$\!\!\!$}{\end{Corollary}}
\newenvironment{remark}{\begin{Remark}$\!\!\!$}{\end{Remark}}

\numberwithin{equation}{section}

\begin{document}

\title{Supersolutions for a class of nonlinear parabolic systems}
%
%
\author{Kazuhiro Ishige\\
Mathematical Institute, Tohoku University\\
Aoba, Sendai 980-8578, Japan,
\\ \\
Tatsuki Kawakami\\
Department of Mathematical Sciences\\
Osaka Prefecture University\\
Sakai 599-8531, Japan
\\ \\
and
\\
\quad
\\ 
Miko{\l}aj Sier\.{z}\c{e}ga\\
Faculty of Mathematics, Informatics and Mechanics\\ 
University of Warsaw\\
Banacha 2, 02-097 Warsaw, Poland\\
}
\date{}
\maketitle
\begin{abstract}
In this paper, 
by using scalar nonlinear parabolic equations, 
we construct supersolutions  
for a class of nonlinear parabolic systems including 
$$
\left\{
\begin{array}{ll}
\partial_t u=\Delta u+v^p,\qquad & x\in\Omega,\,\,\,t>0,\vspace{5pt}\\
\partial_t v=\Delta v+u^q, & x\in\Omega,\,\,\,t>0,\vspace{5pt}\\
u=v=0, & x\in\partial\Omega,\,\,\,t>0,\vspace{5pt}\\
(u(x,0), v(x,0))=(u_0(x),v_0(x)), & x\in\Omega,
\end{array}
\right.
$$
where $p\ge 0$, $q\ge 0$, $\Omega$ is a (possibly unbounded) smooth domain 
in ${\bf R}^N$ 
and both $u_0$ and $v_0$ are nonnegative and locally integrable functions in $\Omega$. 
The supersolutions enable us to obtain optimal 
sufficient conditions for the existence of the solutions 
and optimal lower estimates of blow-up rate of the solutions. 
\end{abstract}
\newpage
\section{Introduction}
This paper is concerned with 
nonnegative solutions of nonlinear parabolic systems.  
Nonlinear parabolic systems have been studied intensively for more than 20 years. 
See \cite{AHV}--\cite{CF}, \cite{DL1}--\cite{FM}, \cite{INS}, \cite{LS}--\cite{QS}, \cite{R, S, Z} and references therein 
(see \cite{QSBook} for a survey).  
Their analysis is generally more complicated than 
that of scalar nonlinear parabolic equations. 
In this paper, by using scalar nonlinear parabolic equations, 
we construct supersolutions of 
\begin{equation}
\label{eq:1.1}
\left\{
\begin{array}{ll}
\partial_t u=\Delta u+v^p,\qquad & x\in\Omega,\,\,\,t>0,\vspace{5pt}\\
\partial_t v=\Delta v+u^q, & x\in\Omega,\,\,\,t>0,\vspace{5pt}\\
u=v=0, & x\in\partial\Omega,\,\,\,t>0,\vspace{5pt}\\
(u(x,0), v(x,0))=(u_0(x),v_0(x)), & x\in\Omega,
\end{array}
\right.
\end{equation}
where $p\ge 0$, $q\ge 0$, $\Omega$ is a (possibly unbounded) smooth domain 
in ${\bf R}^N$ $(N\ge 1)$ 
and both $u_0$ and $v_0$ are nonnegative and locally integrable functions in $\Omega$. 
These supersolutions enable us to obtain optimal sufficient conditions for the existence of 
local-in-time solutions and global-in-time solutions.  
Our arguments are simple and applicable to various nonlinear parabolic systems 
without complicated calculations due to combinations of power and exponential nonlinearities. 
\vspace{5pt}

Problem~\eqref{eq:1.1} is an example of  a simple reaction-diffusion system 
that can be used as a model to describe heat propagation in a two component combustible mixture. 
There are several results on the existence of solutions of \eqref{eq:1.1}. 
Here we recall the following well-known results, 
which were proved in \cite{AHV, EH01, EH02, QS} 
(see also \cite[Section~32]{QSBook}). 
\begin{itemize}
  \item[{\rm (A)}] 
  Let $p,q\ge1$ and $r_1,r_2\in(1,\infty)$. 
  Assume 
  $$
  \max\{P(r_1,r_2),Q(r_1,r_2)\}\le 2,
  $$
  where 
  $$
  P(r_1,r_2):=N\left(\frac{p}{r_2}-\frac{1}{r_1}\right),
  \qquad
  Q(r_1,r_2):=N\left(\frac{q}{r_1}-\frac{1}{r_2}\right).
 $$
  Then, for any $(u_0,v_0)\in L^{r_1}(\Omega)\times L^{r_2}(\Omega)$, 
  problem~\eqref{eq:1.1} possesses a local-in-time solution.   
  \item[{\rm (B)}]  
  Let $0<pq\le 1$. 
  Then problem~\eqref{eq:1.1} possesses a global-in-time solution 
  for any $(u_0,v_0)\in L^\infty(\Omega)\times L^\infty(\Omega)$.
  \item[{\rm (C)}]
  Let $\Omega={\bf R}^N$ and $pq>1$. 
  If 
  \begin{equation}
  \label{eq:1.3}
  \frac{\max\{p,q\}+1}{pq-1}<\frac{N}{2},
  \end{equation}
  then problem~\eqref{eq:1.1} possesses a global-in-time positive solution provided that 
  $(u_0,v_0)\not\equiv(0,0)$ and 
  both $\|u_0\|_{L^{r_1^*}(\Omega)}$ and $\|v_0\|_{L^{r_2^*}(\Omega)}$ are sufficiently small, 
  where
  \begin{equation}
  \label{eq:1.4}
  r_1^*:=\frac{N}{2}\frac{pq-1}{p+1},
  \qquad
  r_2^*:=\frac{N}{2}\frac{pq-1}{q+1}. 
  \end{equation}
  \item[{\rm (D)}]
  Let $\Omega={\bf R}^N$ and $pq>1$.
  If $(u_0,v_0)\not\equiv(0,0)$ and
  $$
  \frac{\max\{p,q\}+1}{pq-1}\ge \frac{N}{2},
  $$
  then problem~\eqref{eq:1.1} admits no global-in-time positive solution.
\end{itemize}
The optimality of assumption~\eqref{eq:1.3} in (C) 
follows from (D). 

In this paper we aim to construct supersolutions 
for a class of nonlinear parabolic systems including \eqref{eq:1.1} 
and to use them to deduce sufficient conditions for the existence of local-in-time solutions. 
In particular, our results for problem~\eqref{eq:1.1} stated in Section~3 improve on (A)--(C). 
Furthermore, as an application of these sufficient conditions, 
we obtain lower estimates on the blow-up rate for the solutions of \eqref{eq:1.1}. 
In Subsection~3.2 and Section~4 we address some possible generalizations to other nonlinear parabolic systems. 

Let us now outline the construction of supersolutions. 
Given $(u,v)$ a positive (classical) solution of \eqref{eq:1.1}, 
we begin by setting $U:=u^\alpha$ and $V:=v^\beta$, where $\alpha\ge 1$ and $\beta\ge 1$. 
Then $(U,V)$ satisfies 
\begin{equation}
\label{eq:1.5}
\left\{
\begin{array}{ll}
\partial_t U=\Delta U+\alpha U^{1-\frac{1}{\alpha}}V^{\frac{p}{\beta}}
-\displaystyle{\frac{\alpha-1}{\alpha}\frac{|\nabla U|^2}{U}}, &x\in\Omega,\,\,t>0,\vspace{7pt}\\
\partial_t V=\Delta V+\beta V^{1-\frac{1}{\beta}}U^{\frac{q}{\alpha}}
-\displaystyle{\frac{\beta-1}{\beta}\frac{|\nabla V|^2}{V}}, &x\in\Omega,\,\,t>0,\vspace{5pt}\\
U=V=0, &x\in\partial\Omega,\,\,t>0,\vspace{3pt}\\
(U(x,0),V(x,0))=(u_0(x)^\alpha,v_0(x)^\beta),\quad &x\in\Omega.
\end{array}
\right.
\end{equation}
Let $(\tilde{U},\tilde{V})$ be a positive solution of 
\begin{equation}
\label{eq:1.6}
\left\{
\begin{array}{ll}
\partial_t \tilde{U}
=\Delta \tilde{U}+\alpha \tilde{U}^{1-\frac{1}{\alpha}}\tilde{V}^{\frac{p}{\beta}},
\qquad & x\in\Omega,\,\,t>0,\vspace{5pt}\\
\partial_t\tilde{V}=\Delta\tilde{V}+\beta\tilde{V}^{1-\frac{1}{\beta}}\tilde{U}^{\frac{q}{\alpha}}, 
& x\in\Omega,\,\,t>0,\vspace{5pt}\\
\tilde{U}=\tilde{V}=0, & x\in\partial\Omega,\,\,t>0,\vspace{5pt}\\
(\tilde{U}(x,0),\tilde{V}(x,0))=(u_0(x)^\alpha,v_0(x)^\beta),\quad &x\in\Omega.
\end{array}
\right.
\end{equation}
Since $\alpha\ge 1$ and $\beta\ge 1$, 
by \eqref{eq:1.5} and \eqref{eq:1.6} 
we see that $(U,V)$ is a subsolution of \eqref{eq:1.6}. 
It follows from the comparison principle that 
$$
\tilde{U}(x,t)\ge u(x,t)^\alpha,
\qquad
\tilde{V}(x,t)\ge v(x,t)^\beta,
\qquad
x\in\Omega,\,\,\,t>0. 
$$
Let $\overline{w}$ be a solution of 
\begin{equation}
\label{eq:1.7}
\left\{
\begin{array}{ll}
\partial_t w=\Delta w+\alpha w^A+\beta w^B,\qquad & x\in\Omega,\,\,t>0,\vspace{5pt}\\
w=0, & x\in\partial\Omega,\,\,t>0,\vspace{5pt}\\
w(x,0)=u_0(x)^\alpha+v_0(x)^\beta,& x\in\Omega,
\end{array}
\right.
\end{equation}
where 
$$
A:=1-\frac{1}{\alpha}+\frac{p}{\beta},
\qquad
B:=1-\frac{1}{\beta}+\frac{q}{\alpha}. 
$$  
Then $(\overline{w},\overline{w})$ is a supersolution of \eqref{eq:1.6}. 
This implies that 
$(\overline{w}^{\frac{1}{\alpha}},\overline{w}^{\frac{1}{\beta}})$ is a supersolution of \eqref{eq:1.1}. 
Therefore, by the comparison principle we obtain 
\begin{equation}
\label{eq:1.8}
0\le u(x,t)^\alpha\le \overline{w}(x,t),
\qquad
0\le v(x,t)^\beta\le \overline{w}(x,t),
\qquad
x\in\Omega,\,\,\,t>0. 
\end{equation}
This supersolution with a suitable choice of $\alpha$ and $\beta$ 
enables us to obtain sufficient conditions 
for the existence of local-in-time solutions and global-in-time solutions 
of problem~\eqref{eq:1.1}. 
Compared with the results in \cite{EH01, EH02, QS}, 
we see that our sufficient conditions are optimal. 
By similar arguments 
we can construct supersolutions of various nonlinear parabolic systems systematically 
and give sufficient conditions for the existence of solutions.  
See Subsection~3.2 and Section~4. 
\vspace{3pt}

The rest of this paper is organized as follows. 
In Section~2 we introduce notation and 
prove some lemmas on the existence of the solutions of \eqref{eq:1.1} 
and a scalar nonlinear parabolic equation. 
In Section~3
we prove two theorems on the existence of the solutions of \eqref{eq:1.1} 
that improve on (A)--(C). 
Furthermore, we give lower estimates on the blow-up rate for the solutions.  
In Section~4 we apply our techniques to parabolic systems with strongly coupled nonlinearities. 
\section{Preliminaries}
We introduce some notation and define the notion of solution for \eqref{eq:1.1}. 
Furthermore, we prove some preliminary lemmas. 
In what follows, $C$ denotes a generic constant. 
\vspace{3pt}

Let $1\le r\le\infty$. 
We define the spaces $L^{r,\infty}(\Omega)$ and $L^{r,\infty}_{{\rm uloc}}(\Omega)$. 
Let $f$ be a measurable function in a smooth domain $\Omega\subset{\bf R}^N$. 
Setting $f=0$ outside of $\Omega$, we define
$$
\mu(\lambda):=\left|\{x\in{\bf R}^N\,:\,|f(x)|>\lambda\}\right|,
\qquad \lambda\ge 0,
$$
i.e. the distribution function of $f$. 
Furthermore, we define 
the non-increasing rearrangement of $f$ by 
$$
f^*(s):=\inf\{\lambda>0\,:\,\mu(\lambda)\le s\}.
$$
Then the spherical rearrangement of $f$ is defined by 
$$
f^\sharp(x):=f^*(c_N|x|^N),
$$
where $c_N$ is the volume of the unit ball in ${\bf R}^N$. 
We define  
$$
L^{r,\infty}(\Omega)
:=\left\{f\,:\, \mbox{$f$ is measurable in $\Omega$},\,\,\, \|f\|_{L^{r,\infty}(\Omega)}<\infty\right\},
$$
where 
$$
\|f\|_{L^{r,\infty}(\Omega)}:=\sup_{s>0}\,s^{1/r}f^{**}(s),
\qquad
f^{**}(s):=\frac{1}{s}\int_0^s f^*(r)\,dr. 
$$
Then the following holds (see e.g., \cite[Section~1.1]{G}).
\begin{itemize}
  \item Let $1<r<\infty$. Then $f\in L^{r,\infty}(\Omega)$ if and only if 
  \begin{equation}
  \label{eq:2.1}
  0\le f^\sharp(x)\le C|x|^{-N/r},\qquad x\in{\bf R}^N.
  \end{equation}
  \item $L^r(\Omega)\subset L^{r,\infty}(\Omega)$ and 
  $L^r(\Omega)\not=L^{r,\infty}(\Omega)$ if $1<r<\infty$ 
  and $L^{r,\infty}(\Omega)=L^r(\Omega)$ if $r\in\{1,\infty\}$. 
  Furthermore, 
  $$
  \|f\|_{L^{r,\infty}(\Omega)}\le\|f\|_{L^r(\Omega)}
  $$
  for any $f\in L^r(\Omega)$, where $1\le r\le\infty$.
  \item Let $1\le r\le\infty$ and let $\{r_j\}_{j=1}^k\subset[1,\infty]$ be such that 
  $$
  \frac{1}{r}=\frac{1}{r_1}+\cdots+\frac{1}{r_k}.
  $$
  Then  
  \begin{equation}
  \label{eq:2.2}
  \biggr\|\prod_{j=1}^kf_j\biggr\|_{L^{r,\infty}(\Omega)}\le C\prod_{j=1}^k\|f_j\|_{L^{r_j,\infty}(\Omega)}
  \end{equation}
  for $f_j\in L^{r_j,\infty}(\Omega)$ and $j=1,2,\dots,k$.
\end{itemize}
For any $x\in{\bf R}^N$ and $R>0$, we put $B(x,R):=\{y\in{\bf R}^N: |x-y|<R\}$.
We define 
\begin{equation*}
\begin{split}
L^{r,\infty}_{{\rm uloc}}(\Omega) & :=\left\{f\in L^r_{{\rm loc}}(\Omega)\,:\,
\sup_{x\in\overline{\Omega}}\,\|f\|_{L^{r,\infty}(\Omega\,\cap\,B(x,1))}<\infty\right\}
\quad\mbox{if}\quad 1\le r<\infty,\\
L^{r,\infty}_{{\rm uloc}}(\Omega) & :=L^\infty(\Omega)\quad\mbox{if}\quad r=\infty.
\end{split}
\end{equation*}
For any $\rho>0$, we set 
$$
|||f|||_{r,\rho}:=
\left\{
\begin{array}{ll}
\displaystyle{\sup_{x\in\overline{\Omega}}}\,\,\|f\|_{L^{r,\infty}(\Omega\,\cap\,B(x,\rho))}
 & \quad\mbox{if}\quad 1\le r<\infty,\vspace{5pt}\\
\|f\|_{L^\infty(\Omega)}
 & \quad\mbox{if}\quad r=\infty,
\end{array}
\right.
$$
which are equivalent norms of $L^{r,\infty}_{{\rm uloc}}(\Omega)$.

We denote by $S(t)$ the Dirichlet heat semigroup on $\Omega$.  
Then, for any $\phi\in L^r(\Omega)$ $(r\ge 1)$, 
$v(t):=S(t)\phi$ represents the unique bounded solution of 
$$
\partial_t v=\Delta v\quad\mbox{in}\quad\Omega\times(0,\infty),
\quad
v=0\quad\mbox{on}\quad\partial\Omega\times(0,\infty),
\quad
v(x,0)=\phi(x)\quad\mbox{in}\quad\Omega. 
$$
We first show the following.
\begin{lemma}
\label{Lemma:2.1}
There exists a positive constant $c_*>0$ such that 
\begin{equation}
\label{eq:2.3}
\|S(t)\varphi\|_{L^\infty(\Omega)}\le c_*t^{-\frac{N}{2r}}|||\varphi|||_{r,\rho},
\qquad
0<t\le\rho^2,
\end{equation}
for any $\varphi\in L^{r,\infty}_{{\rm uloc}}(\Omega)$, 
where $1\le r\le\infty$. 
\end{lemma}
{\bf Proof.}
Let $\{z_j\}\subset{\bf R}^N$ be such that 
\begin{equation}
\label{eq:2.4}
{\bf R}^N=\bigcup_{j=1}^\infty B(z_j,1),
\qquad 
\sum_{j=1}^\infty e^{-\frac{|z_j|^2}{8}}<\infty. 
\end{equation}
Set $u(x,t):=S(t)\varphi$ and $u_\lambda(x,t):=u(\lambda x,\lambda^2 t)$ 
for $\lambda>0$. Let 
$$
\Gamma(x,y,t):=(4\pi t)^{-\frac{N}{2}}\exp\left(-\frac{|x-y|^2}{4t}\right).
$$
It follows from the comparison principle and \eqref{eq:2.4} that 
$$
|u_\lambda(x,1)|\le\int_{{\bf R}^N}\Gamma(x,y,1)|u_\lambda(y,0)|\,dy
\le\sum_{j=1}^\infty\int_{B(x+z_j,1)}\Gamma(x,y,1)|u_\lambda(y,0)|\,dy
$$
for $x\in\lambda^{-1}\Omega$. 
This together with \eqref{eq:2.2} and \eqref{eq:2.4} implies that 
\begin{equation*}
\begin{split}
|u(\lambda x,\lambda^2)|=|u_\lambda(x,1)| 
 & \le\sum_{j=1}^\infty
\|\Gamma(x,\cdot,1)\|_{L^{r',\infty}(B(x+z_j,1))}\|u_\lambda(0)\|_{L^{r,\infty}(B(x+z_j,1)}\\
 & \le C\sum_{j=1}^\infty\|\Gamma(x,\cdot,1)\|_{L^\infty(B(x+z_j,1))}
 |||u_\lambda(0)|||_{r,1}\\
 & \le C\sum_{j=1}^\infty e^{-\frac{|z_j|^2}{8}}|||u_\lambda(0)|||_{r,1}
 \le C |||u_\lambda(0)|||_{r,1},
\end{split}
\end{equation*}
where $1\le r\le\infty$ and $1\le r'\le\infty$ with $1/r+1/r'=1$. 
Then we have 
$$
\|u(\lambda^2)\|_{L^\infty(\Omega)}
\le C|||u_\lambda(0)|||_{r,1}
\le C\lambda^{-\frac{N}{r}}|||\varphi|||_{r,\lambda}. 
$$
Therefore, setting $\lambda=t^{1/2}$, 
we obtain 
$$
\|u(t)\|_{L^\infty(\Omega)}
\le Ct^{-\frac{N}{2r}}|||\varphi|||_{r,t^{1/2}}
\le Ct^{-\frac{N}{2r}}|||\varphi|||_{r,\rho}
$$
for all $t\in(0,\rho^2]$ and $\rho>0$, which implies \eqref{eq:2.3}. 
$\Box$
\vspace{3pt}

Next we define the solution of \eqref{eq:1.1} in $\Omega\times(0,T)$, where $0<T\le\infty$. 
Let $u$ and $v$ be nonnegative measurable functions in $\Omega\times(0,T)$ such that 
$$
u,v\in L^\infty(\tau,T-\tau:L^\infty(\Omega))
$$ 
for all $\tau\in(0,T/2)$.
We say that $(u,v)$ is a subsolution of \eqref{eq:1.1} in $\Omega\times(0,T)$ if 
$(u,v)$ satisfies
\begin{equation*}
\begin{split}
 & u(x,t)\le [S(t)u_0](x)+\int_0^t [S(t-s)v(s)^p](x)\,ds<\infty,\\
 & v(x,t)\le [S(t)v_0](x)+\int_0^t [S(t-s)u(s)^q](x)\,ds<\infty,
\end{split}
\end{equation*}
for almost all $x\in\Omega$ and $t\in(0,T)$. 
Similarly, we say that $(u,v)$ is a supersolution of \eqref{eq:1.1} in $\Omega\times(0,T)$ if 
$(u,v)$ satisfies
\begin{equation*}
\begin{split}
 & \infty>u(x,t)\ge [S(t)u_0](x)+\int_0^t [S(t-s)v(s)^p](x)\,ds,\\
 & \infty>v(x,t)\ge [S(t)v_0](x)+\int_0^t [S(t-s)u(s)^q](x)\,ds,
\end{split}
\end{equation*}
for almost all $x\in\Omega$ and $t\in(0,T)$. 
Furthermore, we say that $(u,v)$ is a solution of \eqref{eq:1.1} in $\Omega\times(0,T)$ 
if  $(u,v)$ is a subsolution and a supersolution of \eqref{eq:1.1} in $\Omega\times(0,T)$. 
Similarly, we define solutions, supersolutions 
and subsolutions of problem~\eqref{eq:1.7}.  
\begin{lemma}
\label{Lemma:2.2}
Let $0<T\le\infty$. Then problem~\eqref{eq:1.1} possesses a solution in $\Omega\times(0,T)$ 
if and only if problem~\eqref{eq:1.1} possesses a supersolution in $\Omega\times(0,T)$. 
This also holds true for problem~\eqref{eq:1.7}. 
\end{lemma}
{\bf Proof.} 
Let $0<T\le\infty$. 
It suffices to prove that 
the existence of a supersolution of \eqref{eq:1.1} implies the existence of a solution of \eqref{eq:1.1}. 
Let $(\tilde{u},\tilde{v})$ be a supersolution of \eqref{eq:1.1} in $\Omega\times(0,T)$. 
Set $u_1=S(t)u_0$ and $v_1(t)=S(t)v_0$. 
Define $(u_n,v_n)$ $(n=2,3,\dots)$ inductively by 
\begin{equation}
\label{eq:2.5}
u_n(t)=S(t)u_0+\int_0^t S(t-s)v_{n-1}(s)^p\,ds,
\quad
v(t)=S(t)v_0+\int_0^t S(t-s)u_{n-1}(s)^q\,ds. 
\end{equation}
Then it follows inductively that 
$$
0\le u_1\le u_2\le\cdots\le u_n\le\cdots\le\tilde{u}<\infty,
\quad
0\le v_1\le v_2\le\cdots\le v_n\le\cdots\le\tilde{v}<\infty,
$$
for almost all $x\in\Omega$ and $t\in(0,T)$. 
This means that 
$$
u(x,t):=\lim_{n\to\infty}u_n(x,t),
\qquad
v(x,t):=\lim_{n\to\infty}v_n(x,t),
$$
exist for almost all $x\in\Omega$ and $t\in(0,T)$. 
Furthermore, by \eqref{eq:2.5} we see that $(u,v)$ is a solution of \eqref{eq:1.1} in $\Omega\times(0,T)$. 
Similarly, we obtain the desired conclusion for problem~\eqref{eq:1.7}. 
Thus Lemma~\ref{Lemma:2.2} follows. 
$\Box$
\vspace{5pt}
\newline
Applying the comparison principle and Lemma~\ref{Lemma:2.2}, 
we obtain the following.
\begin{lemma}
\label{Lemma:2.3}
Assume that there exists a solution of \eqref{eq:1.7} in $\Omega\times(0,T)$ 
for some $\alpha\ge 1$ and $\beta\ge 1$, where $0<T\le\infty$. 
Then problem~\eqref{eq:1.1} possesses a solution in $\Omega\times(0,T)$. 
\end{lemma}

The end of this section
we give sufficient conditions 
for the existence of the solutions of \eqref{eq:1.7} 
by employing the argument in \cite{RS}. 
\begin{lemma}
\label{Lemma:2.4}
Let $A\ge 0$, $B\ge 0$, $\alpha\ge 1$ and $\beta\ge 1$.
Put $w_0=u_0^\alpha+v_0^\beta$.
\begin{itemize}
  \item[{\rm (i)}] Suppose $\max\{A,B\}\le 1$. 
  Then, for any nonnegative function $w_0\in L^1_{{\rm loc}}(\Omega)$ satisfying
  $S(t)w_0\in L^\infty(\Omega)$ for all $t>0$, 
  problem~\eqref{eq:1.7} possesses a global-in-time solution. 
  \item[{\rm (ii)}] Otherwise, let $\max\{A,B\}\ge1$. Then the following holds.
  \begin{itemize}
  \item[{\rm (a)}]
  Let $1\le\sigma\le\max\{A,B\}$ 
  be such that $w_0\in L^\sigma_{{\rm loc}}(\Omega)$ and
  $S(t)w_0^\sigma\in L^\infty(\Omega)$ for all $t>0$. 
  Then there exist positive constants $\gamma_1$ and $T_1$ such that, if 
  \begin{equation}
  \label{eq:2.6}
  \sup_{0<t\le T}
  \biggr[\left\|S(t)w_0^\sigma\right\|_{L^\infty(\Omega)}^{1-\frac{1}{\sigma}}
  \int_0^t\|S(s)w_0^\sigma\|_{L^\infty(\Omega)}^{\frac{\max\{A,B\}}{\sigma}-1}\,ds\biggr]\le\gamma_1
  \end{equation}
  for some $T>0$, then problem~\eqref{eq:1.7} possesses a solution $w$ 
  in $\Omega\times(0,\min\{T,T_1\}]$ such that 
  $$
  0\le w(x,t)\le 2[S(t)w_0^\sigma]^{\frac{1}{\sigma}}+2t
  \quad\mbox{in}\quad\Omega\times(0,\min\{T,T_1\}].
  $$
  \item[{\rm (b)}]
  Let $1\le\sigma\le\min\{A,B\}$ 
  be such that $w_0\in L^\sigma_{{\rm loc}}(\Omega)$ and
  $S(t)w_0^\sigma\in L^\infty(\Omega)$ for all $t>0$. 
  Then there exists a positive constant $\gamma_2$ such that, if 
  $$
  \sup_{0<t\le T}
  \biggr[\left\|S(t)w_0^\sigma\right\|_{L^\infty(\Omega)}^{1-\frac{1}{\sigma}}
  \int_0^t\left\{\|S(s)w_0^\sigma\|_{L^\infty(\Omega)}^{\frac{A}{\sigma}-1}
  +\|S(s)w_0^\sigma\|_{L^\infty(\Omega)}^{\frac{B}{\sigma}-1}\right\}\,ds\biggr]\le\gamma_2
  $$
  for some $T>0$, then 
  problem~\eqref{eq:1.7} possesses a solution $w$ in $\Omega\times(0,T]$ such that 
  $$
  0\le w(x,t)\le 2[S(t)w_0^\sigma]^{\frac{1}{\sigma}}
  \quad\mbox{in}\quad\Omega\times(0,T]. 
  $$
  \end{itemize}
\end{itemize}
\end{lemma}
{\bf Proof.}
We first prove assertion~(ii)-(a). 
Let $\gamma_1$ be a sufficiently small positive constant and assume \eqref{eq:2.6}. 
Set 
\begin{equation}
\label{eq:2.7}
\overline{w}(t):=2[S(t)w_0^\sigma]^{\frac{1}{\sigma}}+2t.
\end{equation}
Since
\begin{equation*}
\begin{split}
F[\overline{w}](t) & :=S(t)w_0+\int_0^tS(t-s)[\alpha \overline{w}(s)^A+\beta \overline{w}(s)^B]\,ds\\
 & \le\left[S(t)w_0^\sigma\right]^{\frac{1}{\sigma}}
 +\int_0^t S(t-s)\left[1+C\overline{w}(s)^{\max\{A,B\}}\right]\,ds\\
 & \le \left[S(t)w_0^\sigma\right]^{\frac{1}{\sigma}}+t
 +C\int_0^t S(t-s)\overline{w}(s)^{\max\{A,B\}}\,ds,\qquad t>0,
\end{split}
\end{equation*}
we have 
\begin{equation*}
\begin{split}
 F[\overline{w}](t)
  & \le\frac{1}{2}\overline{w}(t)+C\int_0^t S(t-s)\overline{w}(s)^{\max\{A,B\}}ds\\
 & \le \frac{1}{2}\overline{w}(t)+C\int_0^t S(t-s)
 \left[[S(s)w_0^\sigma]^{\frac{\max\{A,B\}}{\sigma}}+s^{\max\{A,B\}}\right]\,ds\\
 & \le\frac{1}{2}\overline{w}(t)+C\int_0^t S(t-s)
 \left\|S(s)w_0^\sigma\right\|_{L^\infty(\Omega)}^{\frac{\max\{A,B\}}{\sigma}-1}S(s)w_0^\sigma\,ds
 +Ct^{1+\max\{A,B\}}\\
 & \le\frac{1}{2}\overline{w}(t)
+C\biggr[\int_0^t\left\|S(s)w_0^\sigma\right\|_{L^\infty(\Omega)}^{\frac{\max\{A,B\}}{\sigma}-1}\,ds\biggr]S(t)w_0^\sigma
+Ct^{1+\max\{A,B\}}
\end{split}
\end{equation*} 
for all $t>0$. 
Then it follows from \eqref{eq:2.6} that
$$
F[\overline{w}](t)
\le \frac{1}{2}\overline{w}(t)+C\gamma_1[S(t)w_0^\sigma]^{\frac{1}{\sigma}}+Ct^{1+\max\{A,B\}}
$$
for all $t\in(0,T]$. 
Therefore, taking a sufficiently small $\gamma_1$ if necessary, 
we can find a constant $\tau\in(0,1)$ such that 
$$
F[\overline{w}](t)\le\overline{w},\qquad t\in(0,\min\{T,\tau\}]. 
$$
This implies
that $\overline{w}$ is a supersolution of \eqref{eq:1.7} in $\Omega\times(0,\min\{T,\tau\}]$. 
Then assertion~(ii)-(a) follows from Lemma~\ref{Lemma:2.2}. 
Similarly, setting 
$$
\overline{w}(t):=2[S(t)w_0^\sigma]^{\frac{1}{\sigma}}
$$
instead of \eqref{eq:2.7}, 
we have assertion~(ii)-(b). 

It remans to prove assertion~(i). 
Since $0\le A,B\le 1$, it follows that 
$$
\alpha x^A+\beta x^B\le (\alpha+\beta)(x+1)
$$
for any $x>0$. 
This implies that 
$$
e^{(\alpha+\beta)t}S(t)w_0+e^{(\alpha+\beta)t}-1
$$
is a supersolution of \eqref{eq:1.7} in $\Omega\times(0,\infty)$. 
Then assertion~(i) follows from Lemma~\ref{Lemma:2.2}.  
$\Box$ 
\section{Weakly coupled nonlinear parabolic systems}
In this section, using supersolutions constructed in Section~1 and 
applying Lemmas~\ref{Lemma:2.3} and \ref{Lemma:2.4}, 
we study the existence of local-in-time solutions and global-in-time solutions of \eqref{eq:1.1}. 
Furthermore, we obtain lower estimates on the blow-up rate for the solutions. 
\subsection{Existence of the solutions}
We first give sufficient conditions for the existence of the solutions of \eqref{eq:1.1} 
by using uniformly local $L^{r,\infty}$ spaces. 
In this subsection, we write
$$
P=P(r_1,r_2):=N\left(\frac{p}{r_2}-\frac{1}{r_1}\right)
\qquad\mbox{and}\qquad
Q=Q(r_1,r_2):=N\left(\frac{q}{r_1}-\frac{1}{r_2}\right)
$$
for simplicity.
\begin{theorem}
\label{Theorem:3.1}
Let $u_0$ and $v_0$ be nonnegative measurable functions in $\Omega$ such that 
$u_0\in L^{r_1,\infty}_{{\rm uloc}}(\Omega)$ and $v_0\in L^{r_2,\infty}_{{\rm uloc}}(\Omega)$, 
where $r_1,r_2\in[1,\infty)$. 
\begin{itemize}
  \item[{\rm (i)}] Let $\max\{P,Q\}\le 0$. Then there exists a global-in-time solution of \eqref{eq:1.1}. 
  \item[{\rm (ii)}] Let $W_0:=u_0^{r_1}+v_0^{r_2}$. Assume
  $0<\max\{P,Q\}\le 2$ and $1<r\le\min\{r_1,r_2\}$. 
  Then there exists $\sigma_*>1$ with the following property:  
  for any $1<\sigma\le\sigma_*$, there exist positive constants $\gamma_1$ and $T_1$ such that, if 
  \begin{equation}
  \label{eq:3.1}
  |||W_0|||_{1,\rho}\le\gamma_1\rho^{N\left(1-\frac{2}{\max\{P,Q\}}\right)}
  \end{equation}
  for some $\rho>0$, then problem~\eqref{eq:1.1} possesses a solution~$(u,v)$
  in $\Omega\times(0,\min\{\rho^2,T_1\}]$ satisfying 
  \begin{equation}
  \label{eq:3.2}
  0\le u(x,t)^{r_1}+v(x,t)^{r_2}\le C[S(t)W_0^\frac{\sigma}{r}](x)^{\frac{r}{\sigma}}+C
  \end{equation}
  in $\Omega\times(0,\min\{\rho^2,T_1\}]$. 
  Furthermore, in the case where $P>0$ and $Q>0$, 
  for any $1<\sigma\le\sigma_*$, 
  there exists a positive constant $\gamma_2$ such that, if 
  \begin{equation}
  \label{eq:3.3}
  |||W_0|||_{1,\rho}
  \le\gamma_2\max\left\{\rho^{N\bigr(1-\frac{2}{P}\bigr)},\rho^{N\bigr(1-\frac{2}{Q}\bigr)}\right\}
  \end{equation}
  for some $\rho>0$, then problem~\eqref{eq:1.1} possesses a solution
  in $\Omega\times(0,\rho^2]$ satisfying 
  \begin{equation}
  \label{eq:3.4}
  0\le u(x,t)^{r_1}+v(x,t)^{r_2}\le C[S(t)W_0^\frac{\sigma}{r}](x)^{\frac{r}{\sigma}}
  \end{equation}
 in $\Omega\times(0,\rho^2]$.
  \item[{\rm (iii)}] 
  Let $\max\{P,Q\}>0$. 
  There exist positive constants $\gamma_3$ and $T_2$ such that, if 
  $$
  \int_0^T\|S(s)W_0\|_{L^\infty(\Omega)}^{\frac{\max\{P,Q\}}{N}}\,ds\le\gamma_3
  $$
  for some $T>0$, then problem~\eqref{eq:1.1} possesses a solution
  in $\Omega\times(0,\min\{T,T_2\}]$ satisfying 
  $$
  0\le u(x,t)^{r_1}+v(x,t)^{r_2}\le C[S(t)W_0](x)+C
  $$
  in $\Omega\times(0,\min\{T,T_2\}]$. 
  Furthermore, in the case where $P>0$ and $Q>0$, 
  there exists a positive constant $\gamma_4$ such that, if 
  $$
  \int_0^T\left\{\|S(s)W_0\|_{L^\infty(\Omega)}^{\frac{P}{N}}
  +\|S(s)W_0\|_{L^\infty(\Omega)}^{\frac{Q}{N}}\right\}\,ds\le\gamma_4,
  $$
  then problem~\eqref{eq:1.1} possesses a solution in $\Omega\times(0,T]$ satisfying 
  $$
  0\le u(x,t)^{r_1}+v(x,t)^{r_2}\le C[S(t)W_0](x)
 $$
  in $\Omega\times(0,T]$. 
  \end{itemize}
\end{theorem}
{\bf Proof.} 
Let $1\le r\le\min\{r_1,r_2\}$. Set 
$\alpha:=r_1/r\ge 1$ and $\beta:=r_2/r\ge 1$. 
Then  
\begin{equation}
\label{eq:3.5}
\begin{split}
 & w_0:=u_0^\alpha+v_0^\beta
 =u_0^{\frac{r_1}{r}}+v_0^{\frac{r_2}{r}}\in L^{r,\infty}_{{\rm uloc}}(\Omega),\\
 & A:=1-\frac{1}{\alpha}+\frac{p}{\beta}
=1+r\left(\frac{p}{r_2}-\frac{1}{r_1}\right)
=1+\frac{r}{N}P>0,\\
 & B:=1-\frac{1}{\beta}+\frac{q}{\alpha}
=1+r\left(\frac{q}{r_1}-\frac{1}{r_2}\right)
=1+\frac{r}{N}Q>0. 
\end{split}
\end{equation}
In the case $\max\{P,Q\}\le0$, 
it follows from Lemma~\ref{Lemma:2.4}~(i) that 
problem~\eqref{eq:1.7} possesses a global-in-time positive solution $\overline{w}$. 
Then, by Lemma~\ref{Lemma:2.3}  
we see that problem~\eqref{eq:1.1} possesses a global-in-time positive solution. 
Thus assertion~(i) follows.  

We prove assertion~(ii). Let  
$$
1<r\le\min\{r_1,r_2\},\qquad
1<\sigma\le\max\{A,B\},\qquad
\sigma\le r.
$$  
Since $w_0^\sigma\in L^{\frac{r}{\sigma},\infty}_{{\rm uloc}}(\Omega)$, 
by Lemma~\ref{Lemma:2.1} we have 
$$
\|S(t)w_0^\sigma\|_{L^\infty(\Omega)}\le c_*t^{-\frac{\sigma N}{2r}}
|||w_0^\sigma|||_{\frac{r}{\sigma},\rho}
\le Ct^{-\frac{\sigma N}{2r}}|||W_0|||_{1,\rho}^{\frac{\sigma}{r}},
\qquad 0<t\le\rho^2.
$$
This together with $\sigma>1$, $\max\{P,Q\}\le 2$ and \eqref{eq:3.1} implies that
\begin{equation}
\label{eq:3.6}
\begin{split}
 & \left\|S(t)w_0^\sigma\right\|_{L^\infty(\Omega)}^{1-\frac{1}{\sigma}}
\int_0^t\|S(s)w_0^\sigma\|_{L^\infty(\Omega)}^{\frac{\max\{A,B\}}{\sigma}-1}\,ds\\
 & \le \left[C|||W_0|||_{1,\rho}^{\frac{\sigma}{r}}\right]^{\frac{\max\{A,B\}-1}{\sigma}}
 t^{-\frac{(\sigma-1)N}{2r}}
 \int_0^t s^{-\frac{\max\{P,Q\}}{2}+\frac{(\sigma-1)N}{2r}}\,ds\\
 & \le C|||W_0|||_{1,\rho}^{\frac{\max\{P,Q\}}{N}}
 t^{1-\frac{\max\{P,Q\}}{2}}
 \le C\gamma_1^{\frac{\max\{P,Q\}}{N}},
 \qquad 0<t\le\rho^2.
\end{split}
\end{equation}
Taking a sufficiently small $\gamma_1$ if necessary, 
by Lemma~\ref{Lemma:2.4}~(ii)-(a)
we can find a positive constant $T_1$ such that 
problem~\eqref{eq:1.7} possesses a solution $w$ in $\Omega\times(0,\min\{\rho^2,T_1\}]$ 
such that 
\begin{equation}
\label{eq:3.7}
0<w(x,t)\le 2[S(t)w_0^\sigma]^{\frac{1}{\sigma}}+2t
\end{equation}
in $\Omega\times(0,\min\{\rho^2,T_1\}]$. 
Then, by Lemma~\ref{Lemma:2.3} 
we see that problem~\eqref{eq:1.1} possesses a solution $(u,v)$ in $\Omega\times(0,\min\{\rho^2,T_1\}]$.  
Furthermore, by \eqref{eq:1.8} and \eqref{eq:3.7} we obtain 
\begin{equation}
\label{eq:3.8}
\max\left\{u(x,t)^{\frac{r_1}{r}},v(x,t)^{\frac{r_2}{r}}\right\}
\le w(x,t)\le 2[S(t)w_0^{\sigma}]^{\frac{1}{\sigma}}+2t
\le C[S(t)W_0^\frac{\sigma}{r}]^{\frac{1}{\sigma}}+2t
\end{equation}
in $\Omega\times(0,\min\{\rho^2,T_1\}]$, 
which implies \eqref{eq:3.2}. 

In the case $P>0$ and $Q>0$, 
it follows that $A>1$ and $B>1$. Let $1<\sigma\le\min\{A,B\}$. 
Then, similarly to \eqref{eq:3.6}, by \eqref{eq:3.3} we have 
\begin{equation*}
\begin{split}
 & \left\|S(t)w_0^\sigma\right\|_{L^\infty(\Omega)}^{1-\frac{1}{\sigma}}
\int_0^t\left\{\|S(s)w_0^\sigma\|_{L^\infty(\Omega)}^{\frac{A}{\sigma}-1}
+\|S(s)w_0^\sigma\|_{L^\infty(\Omega)}^{\frac{B}{\sigma}-1}\right\}\,ds\\
 & \le C|||W_0|||_{1,\rho}^{\frac{P}{N}}t^{1-\frac{P}{2}}
+C|||W_0|||_{1,\rho}^{\frac{Q}{N}}t^{1-\frac{Q}{2}}
\le C(\gamma_2^{\frac{P}{N}}+\gamma_2^{\frac{Q}{N}})
\end{split}
\end{equation*}
for all $0<t\le\rho^2$. 
Then, taking a sufficiently small $\gamma_2$ if necessary, 
by Lemma~\ref{Lemma:2.4}~(ii)-(b) 
we can find a solution~$w$ of \eqref{eq:1.7} 
in $\Omega\times(0,\rho^2]$ 
such that 
\begin{equation}
\label{eq:3.9}
0<w(x,t)\le 2[S(t)w_0^\sigma]^{\frac{1}{\sigma}}
\quad\mbox{in}\quad\Omega\times(0,\rho^2]. 
\end{equation}
By Lemma~\ref{Lemma:2.3}  
we can find a solution of \eqref{eq:1.1} in $\Omega\times(0,\rho^2]$. 
Furthermore, similarly to \eqref{eq:3.8}, 
by \eqref{eq:1.8} and \eqref{eq:3.9} we obtain 
$$
\max\left\{u(x,t)^{\frac{r_1}{r}},v(x,t)^{\frac{r_2}{r}}\right\}
\le w(x,t)\le C[S(t)W_0^\frac{\sigma}{r}]^{\frac{1}{\sigma}}
$$
in $\Omega\times(0,\rho^2]$, 
which implies \eqref{eq:3.4}. 
Thus assertion~(ii) follows. 
Assertion~(iii) is also proved by the same argument with $r=\sigma=1$ as in the proof of assertion~(ii) and 
by Lemma~\ref{Lemma:2.4}~(ii). 
Thus Theorem~\ref{Theorem:3.1} follows. 
$\Box$\vspace{5pt}
\newline
It follows that  
$P=Q=2$ if and only if 
\begin{equation*}
r_1=r_1^*:=\frac{N}{2}\frac{pq-1}{p+1}
\quad\mbox{and}\quad 
r_2=r_2^*:=\frac{N}{2}\frac{pq-1}{q+1}. 
\end{equation*}
Then, by Theorem~\ref{Theorem:3.1}~(ii) we have
\begin{corollary}
\label{Corollary:3.1}
Assume that 
$$
r_1^*:=\frac{N}{2}\frac{pq-1}{p+1}>1
\quad\mbox{and}\quad 
r_2^*:=\frac{N}{2}\frac{pq-1}{q+1}>1. 
$$
Then there exists a constant $\delta>0$ such that,  
for any nonnegative measurable functions 
$u_0\in L^{r_1^*,\infty}(\Omega)$ and $v_0\in L^{r_2^*,\infty}(\Omega)$, 
if 
$$
\|u_0\|_{L^{r_1,\infty}(\Omega)}+\|v_0\|_{L^{r_2,\infty}(\Omega)}<\delta,
$$
then problem~\eqref{eq:1.1} possesses a global-in-time solution 
and \eqref{eq:3.4} holds in $\Omega\times(0,\infty)$.
\end{corollary}
Furthermore, by Theorem~\ref{Theorem:3.1} and \eqref{eq:2.1} 
we obtain 
\begin{corollary}
\label{Corollary:3.2}
Let $\Omega$ be a domain in ${\bf R}^N$ such that $0\in\Omega$. 
Let $u_0$ and $v_0$ be nonnegative measurable functions in $\Omega$ such that 
\begin{equation}
\label{eq:3.11}
0\le u_0(x)\le d|x|^{-N/r_1},
\qquad
0\le v_0(x)\le d|x|^{-N/r_2},
\end{equation}
for all $x\in\Omega$, where $d>0$, $r_1>1$ and $r_2>1$. 
In the case $\max\{P,Q\}\le 0$, problem~\eqref{eq:1.1} possesses a global-in-time solution. 
On the other hand, in the case $\max\{P,Q\}>0$, 
there exists a constant $d_*>0$ such that,  
if $0<d\le d_*$, then 
the following holds. 
\begin{itemize}
  \item[{\rm (i)}]
  Let $0<\max\{P,Q\}\le 2$. Then 
  problem~\eqref{eq:1.1} possesses a solution $(u,v)$ in $\Omega\times(0,T)$ for some $T>0$ satisfying 
  $$
  0\le u(x,t)\le C(|x|+t^2)^{-\frac{N}{r_1}}+C,
  \qquad
  0\le v(x,t)\le C(|x|+t^2)^{-\frac{N}{r_2}}+C,
  $$
  in $\Omega\times(0,T)$. 
  \item[{\rm (ii)}] 
  Let $P=Q=2$.
  Then problem~\eqref{eq:1.1} possesses a global-in-time solution $(u,v)$ satisfying 
  $$
  0\le u(x,t)\le C(|x|+t^2)^{-\frac{N}{r_1}},
  \qquad
  0\le v(x,t)\le C(|x|+t^2)^{-\frac{N}{r_2}},
  $$
  in $\Omega\times(0,\infty)$. 
\end{itemize}
\end{corollary}
{\bf Proof.}
Assume \eqref{eq:3.11}. 
It follows from \eqref{eq:2.1} 
that $u_0\in L^{r_1,\infty}(\Omega)$ and $v_0\in L^{r_2,\infty}(\Omega)$. 
Set 
$$
W_0(x):=u_0^{r_1}+v_0^{r_2}. 
$$
For any $1<\sigma\le r\le\min\{r_1,r_2\}$, 
by \eqref{eq:3.11} we have 
$$
0\le W_0(x)^{\frac{\sigma}{r}}\le Cd^\sigma |x|^{-\frac{\sigma N}{r}},
\qquad x\in\Omega,
$$
which implies that 
$$
0\le [S(t)W_0^{\frac{\sigma}{r}}](x)\le Cd^\sigma(|x|+t^2)^{-\frac{\sigma N}{r}},
\qquad x\in\Omega,\,\,t>0. 
$$
Then Corollary~\ref{Corollary:3.2} follows from Theorem~\ref{Theorem:3.1}. 
$\Box$\vspace{7pt}

Next we give a sufficient condition 
for problem~\eqref{eq:1.1} to posses a global-in-time 
positive solution for some initial function. 
Assume the following. 
\begin{equation}
\label{eq:3.12}
\left\{
\begin{array}{l}
 \mbox{There exists a constant $p_*(\Omega)>1$ with the following properties:}\vspace{3pt}\\
 {{\rm(i)}} 
 \quad\,\, \mbox{If $1<\min\{A,B\}\le p_*(\Omega)$, 
  then problem~\eqref{eq:1.7} has no global-in-time}\vspace{3pt}\\ 
 \qquad\,\,\mbox{positive solutions;}\vspace{5pt}\\
  {{\rm(ii)}}
  \quad\mbox{If $p_*(\Omega)<\min\{A,B\}<\infty$, then 
  problem~\eqref{eq:1.7} possesses a global-in-}\vspace{3pt}\\
  \qquad\,\,\mbox{time positive solution for some initial function.} 
\end{array}
\right.
\end{equation}
The critical exponent $p_*(\Omega)$ has been identified  for various domains
{\rm(}see e.g., {\rm\cite{DL, L}}{\rm)}. 
In particular, 
\begin{itemize}
  \item[{\rm (i)}] $p_*(\Omega)=1+2/N$ if $\Omega={\bf R}^N$,
  \item[{\rm (ii)}] $p_*(\Omega)=1+2/(N+1)$ if $\Omega$ is a half space of ${\bf R}^N$,
  \item[{\rm (iii)}] $p_*(\Omega)=1+2/N$ if $\Omega$ is the exterior domain 
  of a compact set in ${\bf R}^N$ and $N\ge 2$. 
\end{itemize} 
(See also {\rm\cite{P}} for {\rm (iii)}.)  
\begin{theorem}
\label{Theorem:3.2}
Let $\Omega$ be a smooth domain in ${\bf R}^N$. 
\begin{itemize}
  \item[{\rm (i)}] If $pq\le 1$, then 
  problem~\eqref{eq:1.1} possesses a global-in-time solution 
  for any initial function $(u_0,v_0)\in L^\infty(\Omega)\times L^\infty(\Omega)$. 
  \item[{\rm (ii)}] Assume \eqref{eq:3.12}. 
  Then problem~\eqref{eq:1.1} possesses a global-in-time positive solution 
  for some initial function if  
  \begin{equation}
  \label{eq:3.13}
  \frac{pq-1}{\max\{p,q\}+1}>p_*(\Omega)-1.
  \end{equation}
\end{itemize}
\end{theorem}
{\bf Proof.}
We prove assertion~(i). Assume $pq\le 1$. 
Let $u_0$, $v_0\in L^\infty(\Omega)$. 
Let $r\ge 1$ be such that $rp\ge 1$. 
Then $u_0\in L^r_{{\rm uloc}}(\Omega)$ and $v_0\in L^{pr}_{{\rm uloc}}(\Omega)$. 
Furthermore, 
$$
P(r,pr)=0,\qquad
Q(r,pr)=\frac{N}{pr}(pq-1)\le 0. 
$$ 
Therefore, assertion~(i) follows from Theorem~{\rm\ref{Theorem:3.1}}~{\rm (i)}. 

We prove assertion~(ii). 
Assume 
$$
A:=1-\frac{1}{\alpha}+\frac{p}{\beta}>p_*(\Omega),
\qquad
B:=1-\frac{1}{\beta}+\frac{q}{\alpha}>p_*(\Omega), 
$$
for some $\alpha\ge 1$ and $\beta\ge 1$. 
This is equivalent to \eqref{eq:3.13}. 
Then assertion~(ii) follows from Lemma~\ref{Lemma:2.3} and the definition of $p_*(\Omega)$. 
Thus the proof is complete. 
$\Box$
\begin{remark}
\label{Remark:3.1}
Assertions~{\rm (B)} and {\rm (C)} follow from Theorem~{\rm\ref{Theorem:3.2}} 
and the fact $p_*({\bf R}^N)=1+2/N$. 
Furthermore, 
assertion~{\rm (A)} follows from Theorem~{\rm\ref{Theorem:3.1}}. 
Indeed, if $u_0\in L^{r_1}(\Omega)$ and $v_0\in L^{r_2}(\Omega)$ 
for some $r_1,r_2\ge 1$, then 
$$
\lim_{\rho\to 0}|||u_0|||_{r_1,\rho}=0
\quad\mbox{and}\quad
\lim_{\rho\to 0}|||u_0|||_{r_2,\rho}=0, 
$$
which means that \eqref{eq:3.1} holds for all sufficiently small $\rho>0$. 
Then assertion~{\rm (A)} follows from Theorem~{\rm\ref{Theorem:3.1}}~{\rm (ii)}. 
\end{remark}
%
\subsection{$k$-component weakly coupled parabolic systems}
Our arguments to problem~\eqref{eq:1.1} 
are applicable to the $k$-component nonlinear parabolic system 
\begin{equation}
\label{eq:3.14}
\left\{
\begin{array}{ll}
\partial_t u_i=\Delta u_i+u_{i+1}^{p_i},\qquad & x\in\Omega,\,\,\,t>0,\vspace{5pt}\\
u_i=0, & x\in\partial\Omega,\,\,\,t>0,\vspace{5pt}\\
u_i(x,0)=u_{i,0}(x), & x\in\Omega,
\end{array}
\right.
\end{equation}
where $k\in\{1,2,3,\dots\}$, 
$i\in\{1,\cdots,k\}$, $p_i\ge 0$ and $u_{k+1}=u_1$. 
Assume that 
$\{u_{i,0}\}$ are nonnegative measurable functions in $\Omega$ 
such that $u_{i,0}\in L^{r_i,\infty}_{\rm uloc}(\Omega)$ $(i=1,\dots,k)$, 
where $r_i\in[1,\infty)$. 
Let $\alpha_1,\alpha_2,\dots,\alpha_k\ge 1$ and $\alpha_{k+1}=\alpha_1$. 
Set $U_i=u_i^{\alpha_i}$ for $i\in\{1,\dots,k\}$ and $U_{k+1}=U_1$. 
Then, similarly to \eqref{eq:1.5}, we have
\begin{equation}
\label{eq:3.15}
\left\{
\begin{array}{ll}
\partial_t U_i=\Delta U_i+\alpha_i U_i^{1-\frac{1}{\alpha_i}}U_{i+1}^{\frac{p_i}{\alpha_{i+1}}}
-\displaystyle{\frac{\alpha_i-1}{\alpha_i}\frac{|\nabla U_i|^2}{U_i}}, &x\in\Omega,\,\,t>0,\vspace{7pt}\\
U_i=0, &x\in\partial\Omega,\,\,t>0,\vspace{3pt}\\
U_i(x,0)=u_{i,0}(x)^{\alpha_i},\quad &x\in\Omega.
\end{array}
\right.
\end{equation}
Set 
$$
A_k:=\max_{1\le i\le k}\left(1-\frac{1}{\alpha_i}+\frac{p_i}{\alpha_{i+1}}\right),
\qquad
B_k:=\min_{1\le i\le k}\left(1-\frac{1}{\alpha_i}+\frac{p_i}{\alpha_{i+1}}\right).
$$
Then there exist positive constants $c_1$ and $c_2$ such that 
$$
\sum_{i=1}^k\alpha_i \xi^{1-\frac{1}{\alpha_i}+\frac{p_i}{\alpha_{i+1}}}
\le c_1\xi^{A_k}+c_2\xi^{B_k},\qquad \xi\in[0,\infty). 
$$
Let $\overline{w}$ be a solution of 
$$
\left\{
\begin{array}{ll}
\partial_t w=\Delta w+c_1 w^{A_k}+c_2 w^{B_k},\qquad & x\in\Omega,\,\,t>0,\vspace{5pt}\\
w=0, & x\in\partial\Omega,\,\,t>0,\vspace{5pt}\\
w(x,0)=\sum_{i=1}^k u_{i,0}(x)^{\alpha_i},& x\in\Omega.
\end{array}
\right.
$$
Then it follows that 
$$
(\overline{w}^{\frac{1}{\alpha_1}},\dots,\overline{w}^{\frac{1}{\alpha_k}})
$$
is a supersolution of \eqref{eq:3.14}. 
Therefore we can apply the arguments in this section to \eqref{eq:3.15} 
with $A$ and $B$ replaced by $A_k$ and $B_k$, respectively. 
In particular, Theorem~\ref{Theorem:3.1} and Corollaries~\ref{Corollary:3.1} and \ref{Corollary:3.2} 
hold with $P(r_1,r_2)$ and $Q(r_1,r_2)$ replaced by 
$$
P_k(r_1,\dots,r_k):=N\max_{1\le i\le k}\left(\frac{p_i}{r_{i+1}}-\frac{1}{r_i}\right),
\quad
Q_k(r_1,\dots,r_k):=N\min_{1\le i\le k}\left(\frac{p_i}{r_{i+1}}-\frac{1}{r_i}\right), 
$$
respectively. 
We leave the details to the reader. 
\subsection{Lower estimates on the blow-up rate}
As an application of the results in subsection~3.1, 
we give lower estimates on the blow-up rate of the solutions of \eqref{eq:1.1} 
by modifying the argument in \cite{IS}, 
which gave lower estimates of the life span of the solutions 
to the hear equation with a nonlinear boundary condition. 
\begin{theorem}
\label{Theorem:3.3} 
Let $\Omega$ be a smooth domain in ${\bf R}^N$ and $pq>1$. 
Let $(u,v)$ be the minimal solution of \eqref{eq:1.1} in $\Omega\times(0,T)$, where $0<T<\infty$, 
such that 
$$
\limsup_{t\to T}\,\left[\|u(t)\|_{L^\infty(\Omega)}+\|v(t)\|_{L^\infty(\Omega)}\right]=\infty. 
$$
Let $r_1^*$ and $r_2^*$ be constants given in \eqref{eq:1.4}
and let $\ell_*\ge 1$ be such that $\ell_*\min\{r_1^*,r_2^*\}\ge 1$. 
Set $\rho(t)=\sqrt{T-t}$. 
Then, for any $r_1\in(\ell_*r_1^*,\infty]$ and $r_2\in(\ell_*r_2^*,\infty]$, 
the following holds.
\begin{itemize}
\item[{\rm(i)}]
There exists a positive constant $c_1$ such that 
$$
\liminf_{t\to T}\,\left\{(T-t)^{\frac{N}{2}(\frac{1}{r_1^*}-\frac{1}{r_1})}|||u(t)|||_{r_1,\rho(t)}
+(T-t)^{\frac{N}{2}(\frac{1}{r_2^*}-\frac{1}{r_2})}|||v(t)|||_{r_2,\rho(t)}\right\}\ge c_1.
$$
\item[{\rm (ii)}]
There exists a positive constant $c_2$ such that
$$
\left\{
\begin{array}{l}
\displaystyle{\limsup_{t\to T}}\,(T-t)^{\frac{p+1}{pq-1}-\frac{N}{2r_1}}|||u(t)|||_{r_1,\rho(t)}\ge c_2,\vspace{3pt}\\
\displaystyle{\limsup_{t\to T}}\,(T-t)^{\frac{q+1}{pq-1}-\frac{N}{2r_2}}|||v(t)|||_{r_2,\rho(t)}\ge c_2.
\end{array}
\right.
$$
\end{itemize}
\end{theorem}
{\bf Proof.} 
We first prove assertion $({\rm i})$. 
Let $\ell>\ell_*$ be such that $r_1\ge r_1':=\ell r_1^*$ and $r_2\ge r_2':=\ell r_2^*$. 
Then 
$$
P(r_1',r_2')=Q(r_1',r_2')=N\left(\frac{p}{\ell r_2^*}-\frac{1}{\ell r_1^*}\right)=\frac{2}{\ell}.
$$
Since the minimal solution $(u,v)$ blows up at $t=T$, 
by Theorem~\ref{Theorem:3.1}~(ii) 
we can find constants $\gamma_1>0$ and $\delta\in(0,T)$ such that 
\begin{equation}
\label{eq:3.16}
|||u(t)^{r_1'}|||_{1,\rho(t)}+|||v(t)^{r_2'}|||_{1,\rho(t)}>\gamma_1\rho(t)^{-N(\ell-1)}
\end{equation}
for all $t\in(T-\delta,T)$. 
Indeed, if not, then we deduce from Theorem~\ref{Theorem:3.1}~(ii) 
that problem~\eqref{eq:1.1} possesses a solution of \eqref{eq:1.1} not blowing up at $t=T$. 

On the other hand, 
it follows from \eqref{eq:2.2} that 
\begin{equation}
\label{eq:3.17}
\begin{split}
|||u(t)^{r_1'}|||_{1,\rho(t)} 
 & =\sup_{x\in\overline{\Omega}}\|u(t)\|_{L^{r_1',\infty}(B(x,\rho(t)))}^{r_1'}\\
 & \le C(T-t)^{\frac{N}{2}\left(1-\frac{r_1'}{r_1}\right)}
 \,\sup_{x\in\Omega}\|u(t)\|_{L^{r_1,\infty}(B(x,\rho(t)))}^{r_1'}\\
 & =C(T-t)^{\frac{N}{2}\left(1-\frac{\ell r_1^*}{r_1}\right)}|||u(t)|||_{r_1,\rho(t)}^{\ell r_1^*}
\end{split}
\end{equation}
and 
\begin{equation}
\label{eq:3.18}
|||v(t)^{r_2'}|||_{1,\rho(t)}
\le C(T-t)^{\frac{N}{2}\left(1-\frac{\ell r_2^*}{r_2}\right)}|||v(t)|||_{r_2,\rho(t)}^{\ell r_2^*}
\end{equation}
for all $t\in(0,T)$. 
By \eqref{eq:3.16}, \eqref{eq:3.17} and \eqref{eq:3.18} we obtain 
\begin{equation}
\label{eq:3.19}
(T-t)^{-\frac{N\ell r_1^*}{2r_1}}|||u(t)|||_{r_1,\rho(t)}^{\ell r_1^*}
+(T-t)^{-\frac{N\ell r_2^*}{2r_2}}|||v(t)|||_{r_2,\rho(t)}^{\ell r_2^*}
\ge\gamma_2\rho(t)^{-N\ell}
\end{equation}
for $t\in(T-\delta,T)$, where $\gamma_2$ is a positive constant. 

Let $t\in(T-\delta,T)$ and assume that 
$$
(T-t)^{-\frac{N\ell r_1^*}{2r_1}}|||u(t)|||_{r_1,\rho(t)}^{\ell r_1^*}
\ge(T-t)^{-\frac{N\ell r_2^*}{2r_2}}|||v(t)|||_{r_2,\rho(t)}^{\ell r_2^*}.
$$
Then, by \eqref{eq:3.19} we have 
$$
(T-t)^{-\frac{N\ell r_1^*}{2r_1}}|||u(t)|||_{r_1,\rho(t)}^{\ell r_1^*}
\ge\frac{\gamma_2}{2}(T-t)^{-\frac{N\ell}{2}},
$$
which implies that 
\begin{equation}
\label{eq:3.20}
|||u(t)|||_{r_1,\rho(t)}\ge\left(\frac{\gamma_2}{2}\right)^{\frac{1}{\ell r_1^*}}(T-t)^{-\frac{N}{2r_1^*}+\frac{N}{2r_1}}
=\left(\frac{\gamma_2}{2}\right)^{\frac{1}{\ell r_1^*}}(T-t)^{-\frac{p+1}{pq-1}+\frac{N}{2r_1}}. 
\end{equation}
If not, we obtain 
\begin{equation}
\label{eq:3.21}
|||v(t)|||_{r_2,\rho(t)}\ge
\left(\frac{\gamma_2}{2}\right)^{\frac{1}{\ell r_2^*}}(T-t)^{-\frac{q+1}{pq-1}+\frac{N}{2r_2}}. 
\end{equation}
We deduce from \eqref{eq:3.20} and \eqref{eq:3.21} that 
\begin{equation*}
\begin{split}
 & (T-t)^{\frac{p+1}{pq-1}-\frac{N}{2r_1}}|||u(t)|||_{r_1,\rho(t)}
+(T-t)^{\frac{q+1}{pq-1}-\frac{N}{2r_2}}|||v(t)|||_{r_2,\rho(t)}\\
 & \ge\min\biggr\{\left(\frac{\gamma_2}{2}\right)^{\frac{1}{\ell r_1^*}},\left(\frac{\gamma_2}{2}\right)^{\frac{1}{\ell r_2^*}}\biggr\}>0
\end{split}
\end{equation*}
for all $t\in(T-\delta,T)$,
which implies assertion $({\rm i})$. 

We prove assertion $(\rm ii)$ by contradiction. 
Let $\epsilon$ be a sufficiently small positive constant. 
Let $r_1>\ell_* r_1^*$ and assume 
\begin{equation}
\label{eq:3.22}
\sup_{t\in(T-\delta,T)}(T-t)^{\frac{p+1}{pq-1}-\frac{N}{2r_1}}|||u(t)|||_{r_1,\rho(t)}\le\epsilon. 
\end{equation}
It follows from \eqref{eq:1.4} that 
\begin{equation}
\label{eq:3.23}
\frac{p+1}{pq-1}-\frac{N}{2r_1}>0.
\end{equation}
Let $r_2':=r_2^*r_1/r_1^*$. Then $r_1>1$, $r_2'>\ell_*r_2^*\ge 1$ and 
\begin{equation}
\label{eq:3.24}
\frac{N}{2}\left(\frac{q}{r_1}-\frac{1}{r_2'}\right)=\frac{Nr_1^*}{2r_1}\left(\frac{q}{r_1^*}-\frac{1}{r_2^*}\right)
=\frac{r_1^*}{r_1}\in(0,1). 
\end{equation}
Put $\tau:=T-\delta$ and $\eta:=(T-t)/2$. 
Then, by Lemma~\ref{Lemma:2.1} and \eqref{eq:3.22}--\eqref{eq:3.24} we obtain 
\begin{equation*}
\begin{split}
 & |||v(t)|||_{r_2',\rho(t)} \\
 & \le|||S(t-\tau)v(\tau)|||_{r_2',\rho(t)}
+\int_\tau^t|||S(t-s)u(s)^q|||_{r_2,\rho(t)}\,ds\\
 & \le C|||v(\tau)|||_{r_2',(T-t)^{\frac{1}{2}}}+C\int_\tau^t(t-s)^{-\frac{N}{2}\left(\frac{q}{r_1}-\frac{1}{r_2'}\right)}|||u(s)|||_{r_1,\rho(t)}^q\,ds\\
 & \le C|||v(\tau)|||_{r_2',T^\frac{1}{2}}+C\epsilon^q\left\{\int_\tau^{t-\eta}+\int_{t-\eta}^t\right\}(t-s)^{-\frac{N}{2}\left(\frac{q}{r_1}-\frac{1}{r_2'}\right)}
 (T-s)^{-\frac{pq+q}{pq-1}+\frac{Nq}{2r_1}}\,ds\\
 & \le C|||v(\tau)|||_{r_2',T^\frac{1}{2}}
 +C\epsilon^q\biggr[
 \eta^{-\frac{N}{2}\left(\frac{q}{r_1}-\frac{1}{r_2'}\right)}(T-t+\eta)^{-\frac{pq+q}{pq-1}}\int_\tau^{t-\eta}(T-s)^\frac{Nq}{2r_1}\,ds\\
 & \qquad\qquad\qquad\qquad\qquad\qquad\qquad
+(T-t)^{-\frac{pq+q}{pq-1}+\frac{Nq}{2r_1}}\int_{t-\eta}^t(t-s)^{-\frac{N}{2}\left(\frac{q}{r_1}-\frac{1}{r_2'}\right)}\,ds\biggr]\\
 & \le C|||v(\tau)|||_{r_2',T^\frac{1}{2}}+C\epsilon^q(T-t)^{-\frac{q+1}{pq-1}+\frac{N}{2r_2'}}
\end{split}
\end{equation*}
for all $0<\tau<t<T$. This implies that 
\begin{equation}
\label{eq:3.25}
\limsup_{t\to T}\,(T-t)^{\frac{q+1}{pq-1}-\frac{N}{2r_2'}}|||v(t)|||_{r_2',\rho(t)}\le\epsilon^q.
\end{equation}
Taking a sufficiently small $\epsilon>0$ if necessary, 
we deduce from \eqref{eq:3.22} and \eqref{eq:3.25} that 
\begin{equation*}
\begin{split}
 & \limsup_{t\to T}\,
\left[(T-t)^{\frac{p+1}{pq-1}-\frac{N}{2r_1}}|||u(t)|||_{r_1,\rho(t)}
+(T-t)^{\frac{q+1}{pq-1}-\frac{N}{2r_2'}}|||v(t)|||_{r_2',\rho(t)}\right]\\
 & \le\epsilon+C\epsilon^q<c_1,
\end{split}
\end{equation*}
where $c_1$ is the constant given in assertion~(i). 
This contradicts assertion~(i), which means  that 
$$
\limsup_{t\to T}\,(T-t)^{\frac{p+1}{pq-1}-\frac{N}{2r_1}}\|u(t)\|_{r_1,\rho(t)}\ge c
$$
for some $c>0$. 
Similarly, we have 
$$
\limsup_{t\to T}\,(T-t)^{\frac{q+1}{pq-1}-\frac{N}{2r_2}}\|v(t)\|_{r_2,\rho(t)}\ge c'
$$
for some $c'>0$. 
Thus assertion~(ii) follows and the proof of Theorem~\ref{Theorem:3.3} is complete. 
$\Box$\vspace{5pt}
\newline
It follows from Theorem~\ref{Theorem:3.3}~(ii) that 
\begin{equation}
\label{eq:3.26}
\limsup_{t\to T}\,(T-t)^{\frac{p+1}{pq-1}}\|u(t)\|_{L^\infty(\Omega)}>0,
\quad
\limsup_{t\to T}\,(T-t)^{\frac{q+1}{pq-1}}\|v(t)\|_{L^\infty(\Omega)}>0.
\end{equation}
(See also \cite[Proposition~3.4]{Z}.)  
For upper estimates of blow-up rate of the solutions of \eqref{eq:1.1}, 
see e.g., \cite{AHV, CM, CF, DL2, QSBook}, which show that 
\eqref{eq:3.26} gives the optimal lower estimate on the blow-up rate of the solutions. 
\section{Parabolic systems with strongly coupled nonlinearities}
In this section we apply the methods of Sections~2 and 3 to  
obtain sufficient conditions for the existence of the solutions 
for parabolic systems with strongly coupled nonlinearities. 
\subsection{Strongly coupled power nonlinearities}
Consider the parabolic system with strongly coupled power nonlineaities  
\begin{equation}
\label{eq:4.1}
\left\{
\begin{array}{ll}
\partial_t u=\Delta u+u^{p_1}v^{p_2},\qquad & x\in\Omega,\,\,\,t>0,\vspace{5pt}\\
\partial_t v=\Delta v+u^{q_1}v^{q_2}, & x\in\Omega,\,\,\,t>0,\vspace{5pt}\\
u=v=0, & x\in\partial\Omega,\,\,\,t>0,\vspace{5pt}\\
(u(x,0), v(x,0))=(u_0(x),v_0(x)), & x\in\Omega,
\end{array}
\right.
\end{equation}
where $p_i\ge 0$, $q_i\ge 0$ $(i=1,2)$, 
$\Omega$ is a (possibly unbounded) smooth domain in ${\bf R}^N$ $(N\ge 1)$ 
and both $u_0$ and $v_0$ are nonnegative and locally integrable functions in $\Omega$. 
Problem~\eqref{eq:4.1} was studied in \cite{C, DL1, EL, FM, LS}, 
however, compared with problems~\eqref{eq:1.1} and \eqref{eq:1.7}, 
much less is known about the conditions for the existence of the solutions of \eqref{eq:4.1}. 
In this subsection we apply the arguments in Section~3 to problem~\eqref{eq:4.1} 
and obtain sufficient conditions for the existence of the solutions.  
Furthermore, we give lower estimates on the blow-up rate of the solutions. 

Let $(u,v)$ be a positive (classical) solution of problem~\eqref{eq:4.1}. 
Set $U:=u^\alpha$ and $V:=v^\beta$ for some $\alpha\ge 1$ and $\beta\ge 1$. 
Then $(U,V)$ satisfies 
$$
\left\{
\begin{array}{ll}
\partial_t U=\Delta U+\alpha U^{1+\frac{p_1-1}{\alpha}}V^{\frac{p_2}{\beta}}
-\displaystyle{\frac{\alpha-1}{\alpha}\frac{|\nabla U|^2}{U}}, &x\in\Omega,\,\,t>0,\vspace{7pt}\\
\partial_t V=\Delta V+\beta U^{\frac{q_1}{\alpha}}V^{1+\frac{q_2-1}{\beta}}
-\displaystyle{\frac{\beta-1}{\beta}\frac{|\nabla V|^2}{V}}, &x\in\Omega,\,\,t>0,\vspace{5pt}\\
U=V=0, &x\in\partial\Omega,\,\,t>0,\vspace{3pt}\\
(U(x,0),V(x,0))=(u_0(x)^\alpha,v_0(x)^\beta),\quad &x\in\Omega.
\end{array}
\right.
$$
Put
$$
A=1+\frac{p_1-1}{\alpha}+\frac{p_2}{\beta},
\qquad
B=1+\frac{q_1}{\alpha}+\frac{q_2-1}{\beta}. 
$$
Let $\overline{w}$ be a solution of \eqref{eq:1.7}. 
Similarly to \eqref{eq:1.1}, we immediately see that
$(\overline{w}^{\frac{1}{\alpha}},\overline{w}^{\frac{1}{\beta}})$
is a supersolution of \eqref{eq:4.1}. 
Then we apply the same arguments as in Sections~2 and 3 
to obtain the following theorems.
\begin{theorem}
\label{Theorem:4.1}
Let $\Omega$ be a smooth domain in ${\bf R}^N$. 
Consider problem~\eqref{eq:4.1}.
Then the same statements as in Theorem~{\rm\ref{Theorem:3.1}}  
and Corollaries~{\rm\ref{Corollary:3.1}} and {\rm\ref{Corollary:3.2}} hold with $P$ and $Q$ 
replaced by 
$$
\tilde{P}:=N\left(\frac{p_1-1}{r_1}+\frac{p_2}{r_2}\right),\qquad
\tilde{Q}:=N\left(\frac{q_1}{r_1}+\frac{q_2-1}{r_2}\right),
$$
respectively. 
\end{theorem}
{\bf Proof.}
Let $1\le r\le\min\{r_1,r_2\}$, $\alpha=r_1/r$ and $\beta=r_2/r$. 
Then, similarly to \eqref{eq:3.5}, we have 
$$
A=1+\frac{p_1-1}{\alpha}+\frac{p_2}{\beta}=1+\frac{r}{N}\tilde{P},
\qquad
B=1+\frac{q_1}{\alpha}+\frac{q_2-1}{\beta}=1+\frac{r}{N}\tilde{Q}.
$$
Then Theorem~\ref{Theorem:4.1} follows from the same arguments in Section~3. 
$\Box$\vspace{5pt}
\newline
We remark that  
$\tilde{P}=\tilde{Q}=2$ if and only if 
$$
r_1=\tilde{r}_1^*:=\frac{N}{2}\frac{\delta}{1-q_2+p_2},
\qquad
r_2=\tilde{r}_2^*:=\frac{N}{2}\frac{\delta}{1-p_1+q_1}, 
$$
where $\delta:=q_1p_2-(p_1-1)(q_2-1)$. 
On the other hand, 
by Lemmas~\ref{Lemma:2.3} and \ref{Lemma:2.4} 
we immediately obtain the following theorem.
\begin{theorem}
\label{Theorem:4.2}
Let $\Omega$ be a smooth domain in ${\bf R}^N$. 
\begin{itemize}
\item[{\rm (i)}] Assume that 
\begin{equation}
\label{eq:4.2}
\frac{p_1-1}{\alpha}+\frac{p_2}{\beta}\le 0,
\qquad
\frac{q_2-1}{\beta}+\frac{q_1}{\alpha}\le 0,
\end{equation}
for some $\alpha\ge 1$ and $\beta\ge 1$. 
Then problem~\eqref{eq:4.1} possesses a global-in-time solution 
for any initial function $(u_0,v_0)\in L^\infty(\Omega)\times L^\infty(\Omega)$. 
\item[{\rm (ii)}] 
Assume \eqref{eq:3.12}. Furthermore, assume that 
\begin{equation}
\label{eq:4.3}
1+\frac{p_1-1}{\alpha}+\frac{p_2}{\beta}>p_*(\Omega),
\qquad
1+\frac{q_2-1}{\beta}+\frac{q_1}{\alpha}>p_*(\Omega),
\end{equation}
for some $\alpha\ge 1$ and $\beta\ge 1$. 
Then problem~\eqref{eq:4.1} possesses a global-in-time positive solution 
for some initial function $(u_0,v_0)$. 
\end{itemize}
\end{theorem}
\begin{remark}
\label{Remark:4.1} 
\eqref{eq:4.2} holds if and only if 
$$
p_1\le 1,
\quad
q_2\le 1,
\quad
\delta\le 0. 
$$
This is the same condition as in Theorem~{\rm 3.2} in {\rm\cite{C}} 
and Theorem~{\rm 5}~{\rm II}-{\rm A} in {\rm\cite{EL}}.  
Furthermore, \eqref{eq:4.3} holds if and only if 
one of the following holds: 
\begin{itemize}
  \item $p_1+p_2>p_*(\Omega)$ and  $q_1+q_2>p_*(\Omega)$;
  \item $p_1<1$, $p_1+p_2\le p_*(\Omega)<q_1+q_2$ and $\delta>(p_*(\Omega)-1)(1-p_1+q_1)$;
  \item $q_2<1$, $q_1+q_2\le p_*(\Omega)<p_1+p_2$ and $\delta>(p_*(\Omega)-1)(1-q_2+p_2)$. 
\end{itemize}
In the case $\Omega={\bf R}^N$, 
this is the same conditions as in Theorem~{\rm 5}~{\rm I}-{\rm A} and {\rm II}-{\rm B} in {\rm\cite{EL}} 
and it is the optimal condition for the existence of global-in-time positive solutions of \eqref{eq:4.1}. 
\end{remark}
On the other hand, 
similarly to Theorem~\ref{Theorem:3.3}~(i), 
we have:
\begin{theorem}
\label{Theorem:4.3} 
Assume $\tilde{r}_1^*>0$ and $\tilde{r}_2^*>0$ .
Let $(u,v)$ be the minimal solution of \eqref{eq:4.1} in $\Omega\times(0,T)$, where $0<T<\infty$, 
such that 
$$
\limsup_{t\to T}\,\left[\|u(t)\|_{L^\infty(\Omega)}+\|v(t)\|_{L^\infty(\Omega)}\right]=\infty. 
$$
Let $\ell_*\ge 1$ be such that $\ell_*\min\{\tilde{r}_1^*,\tilde{r}_2^*\}\ge 1$. 
Then, for any $r_1\in(\ell_*r_1^*,\infty]$ and $r_2\in(\ell_*r_2^*,\infty]$, 
there exists a positive constant $c_1$ such that 
$$
\liminf_{t\to T}\,\left\{(T-t)^{\frac{N}{2}(\frac{1}{r_1^*}-\frac{1}{r_1})}|||u(t)|||_{r_1,\rho(t)}
+(T-t)^{\frac{N}{2}(\frac{1}{r_2^*}-\frac{1}{r_2})}|||v(t)|||_{r_2,\rho(t)}\right\}\ge c_1,
$$
where $\rho(t)=\sqrt{T-t}$. In particular, 
\begin{equation}
\label{eq:4.4}
\liminf_{t\to T}\,\left\{(T-t)^{\frac{1-q_2+p_2}{\delta}}||u(t)||_{L^\infty(\Omega)}
+(T-t)^{\frac{1-p_1+q_1}{\delta}}||v(t)||_{L^\infty(\Omega)}\right\}\ge c_1. 
\end{equation}
\end{theorem}
For upper estimates on the blow-up rate of the solutions of \eqref{eq:4.1}, 
see e.g., \cite[Theorem~1.1]{R}, which shows that 
the lower estimate on the blow-up rate \eqref{eq:4.4} is optimal. 
\begin{remark}
\label{Remark:4.2}
Similarly to Subsection~{\rm 3.2}, 
we can apply the arguments in this subsection to 
the $k$-component nonlinear parabolic system
$$
\left\{
\begin{array}{ll}
{\displaystyle \partial_t u_i=\Delta u_i+\prod_{j=1}^ku_{j}^{p_{i,j}}},\qquad & x\in\Omega,\,\,\,t>0,\vspace{5pt}\\
u_i=0, & x\in\partial\Omega,\,\,\,t>0,\vspace{5pt}\\
u_i(x,0)=u_{i,0}(x), & x\in\Omega,
\end{array}
\right.
$$
where $i\in\{1,\cdots,k\}$, $p_{i,j}\ge 0$ $(j=1,\cdots,k)$ 
and $\{u_{i,0}\}$ are nonnegative and locally integrable functions in $\Omega$. 
We leave the details to the reader. 
\end{remark}
\subsection{Strongly coupled exponential nonlinearities}
Consider the parabolic system with strongly coupled exponential nonlinearities
\begin{equation}
\label{eq:4.5}
\left\{
\begin{array}{ll}
\partial_t u=\Delta u+e^{p_1 u}e^{p_2v},\qquad & x\in{\bf R}^N,\,\,\,t>0,\vspace{5pt}\\
\partial_t v=\Delta v+e^{q_1 u}e^{q_2 v}, & x\in{\bf R}^N,\,\,\,t>0,\vspace{5pt}\\
(u(x,0), v(x,0))=(u_0(x),v_0(x)),\qquad & x\in{\bf R}^N,
\end{array}
\right.
\end{equation}
where $p_i\ge 0$, $q_i\ge 0$ $(i=1,2)$ 
and both $u_0$ and $v_0$ are locally integrable functions in ${\bf R}^N$. 
Set $\tilde{u}:=e^u$ and $\tilde{v}=e^v$. Then $(\tilde{u},\tilde{v})$ satisfies 
$$
\left\{
\begin{array}{ll}
\partial_t\tilde{u}=\Delta\tilde{u}+\tilde{u}^{p_1+1}\tilde{v}^{p_2}
-\displaystyle{\frac{|\nabla \tilde{u}|^2}{\tilde{u}}},
\qquad & x\in{\bf R}^N,\,\,\,t>0,\vspace{8pt}\\
\partial_t\tilde{v}=\Delta\tilde{v}+\tilde{u}^{q_1}\tilde{v}^{q_2+1}
-\displaystyle{\frac{|\nabla \tilde{v}|^2}{\tilde{v}}}, & x\in{\bf R}^N,\,\,\,t>0,\vspace{8pt}\\
(\tilde{u}(x,0), \tilde{v}(x,0))=(e^{u_0(x)},e^{v_0(x)}),\qquad & x\in{\bf R}^N.
\end{array}
\right.
$$
Similarly to Lemma~\ref{Lemma:2.2}, 
we see that problem~\eqref{eq:4.5} possesses a solution $(u,v)$ in ${\bf R}^N\times(0,T)$, 
where $0<T\le\infty$, if the problem 
$$
\left\{
\begin{array}{ll}
\partial_t\hat{u}=\Delta\hat{u}+\hat{u}^{p_1+1}\tilde{v}^{p_2},
\qquad & x\in{\bf R}^N,\,\,\,t>0,\vspace{5pt}\\
\partial_t\hat{v}=\Delta\hat{v}+\hat{u}^{q_1}\hat{v}^{q_2+1},
 & x\in{\bf R}^N,\,\,\,t>0,\vspace{5pt}\\
(\hat{u}(x,0), \hat{v}(x,0))=(e^{u_0(x)},e^{v_0(x)}),\qquad & x\in{\bf R}^N
\end{array}
\right.
$$
possesses a solution $(\hat{u},\hat{v})$ in ${\bf R}^N\times(0,T)$. 
Then we can apply the arguments in Subsection~4.1 
and obtain sufficient conditions for the existence of the solutions of \eqref{eq:4.5}. 
We leave the details to the reader again.  
\,\\

\noindent
{\bf Acknowledgements.}
The first author was supported
by the Grant-in-Aid for for Scientific Research (A)(No.~15H02058),
from Japan Society for the Promotion of Science.
The second author was supported by the Grant-in-Aid for Young Scientists (B)
(No.~24740107)
from Japan Society for the Promotion of Science
and by the JSPS Program for Advancing Strategic International Networks 
to Accelerate the Circulation of Talented Researchers 
``Mathematical Science of Symmetry, Topology and Moduli, 
Evolution of International Research Network based on OCAMI". 
The third author was partially supported by WCMCS, Warsaw.  

\end{document}